\documentclass[english]{article}
\usepackage[T1]{fontenc}
\usepackage[latin9]{inputenc}
\usepackage{geometry}
\usepackage{amstext}
\usepackage{amssymb}
\usepackage{esint}
\usepackage{babel}
\begin{document}
\title{Proof of the Riemann Hypothesis}
\author{Björn Tegetmeyer}
\date{11.10.2023}
\maketitle

\section*{Abstract}

The Riemann hypothesis, stating that the real part of all non-trivial
zero points of the zeta function must be $\frac{1}{2}$, is one of
the most important unproven hypotheses in number theory. In this paper
we will prove the Riemann hypothesis by using the integral representation
$\zeta(s)=\frac{s}{s-1}-s\int_{1}^{\infty}\frac{x-\lfloor x\rfloor}{x^{s+1}}\,\text{d}x$
and solving the integral for the real- and imaginary part of the zeta
function. 

\section{Introduction}

In 1859 Bernhard Riemann found one of the most eminent mathematical
problems of our time: In his paper ``On the Number of Primes Less
Than a Given Magnitude'' [4] he published the assumption
that all non-trivial zero-points of the zeta function extended to
the range of complex numbers $\mathbb{C}$ have a real part of $\frac{1}{2}$,
noting the demand of a strict proof for this. Ever since David Hilbert
in 1900 added this problem to his list of the 23 most important problems
of 20$^{\text{th}}$ century, mathematicians have been working on
finding evidence for Riemanns hypothesis. This paper aims to provide
the proof and fill this gap in modern mathematics. 

\section{Proof of the Riemann Hypothesis}

The zeta-function $\zeta(s)$ in the complex range $s\in\mathbb{C}$
for a positive real-part of $s$ can be formulated as integral representation
\begin{equation}
\zeta(s)=\frac{s}{s-1}-s\int_{1}^{\infty}\frac{x-\lfloor x\rfloor}{x^{s+1}}\,\text{d}x\label{eq:zeta_s}
\end{equation}
 with $s\in\mathbb{C}$, where $s$ can be expressed by $s=a+ib;\,a,b\in\mathbb{R}$
and $0<a<1$ as well as $0<b$. Be $s_{0}$ a zero point of the zeta
function. From [1] we know, that the zeta-function
is symmetrical in a way that $\zeta(s_{0})=0\Leftrightarrow\zeta(1-s_{0})=0$
for all zero-points $s_{0}\in\mathbb{C}$ (see appendix for details).
In accordance to equation \ref{eq:zeta_s} we can write $\zeta(1-s)$
as 
\begin{equation}
\zeta(1-s)=\frac{1-s}{-s}-(1-s)\int_{1}^{\infty}\frac{x-\lfloor x\rfloor}{x^{2-s}}\,\text{d}x\label{eq:zeta_1-s}
\end{equation}
 The Riemann hypothesis states, that the real part of $s_{0}$ would
be $\frac{1}{2}$ for all non-trivial zero-points of zeta (i.e. all
zero points of zeta with a positive real part). Furthermore, from
[2] we know, that $0<\Re(s_{0})<1$. Inserting $s_{0}=a+ib$
into $\zeta(s)$, using $x^{-1-a-ib}=x^{-1-a}\cos\left(b\ln(x)\right)-ix^{-1-a}\sin\left(b\ln(x)\right)$
and defining $\{x\}:=x-\lfloor x\rfloor$ we get the following two
equations \ref{eq:3} and \ref{eq:4} out off $\zeta(s)$ and $\zeta(1-s)$:
\begin{equation}
\begin{array}{rcl}
\frac{1}{a+ib-1} & = & \int_{1}^{\infty}\frac{\{x\}}{x^{a+ib+1}}\,\text{d}x\\
\Leftrightarrow\frac{a-1}{(a-1)^{2}+b^{2}}-i\frac{b}{(a-1)^{2}+b^{2}} & = & \int_{1}^{\infty}\frac{\{x\}}{x^{1+a}}\left(\cos(b\ln(x))-i\sin(b\ln(x))\right)\,\text{d}x
\end{array}\label{eq:3}
\end{equation}
 and 
\begin{equation}
\begin{array}{rcl}
\frac{1}{-a-ib} & = & \int_{1}^{\infty}\frac{\{x\}}{x^{2-a-ib}}\,\text{d}x\\
\Leftrightarrow\frac{-a}{a^{2}+b^{2}}+i\frac{b}{a^{2}+b^{2}} & = & \int_{1}^{\infty}\frac{\{x\}}{x^{2-a}}\left(\cos(b\ln(x))+i\sin(b\ln(x))\right)\,\text{d}x
\end{array}\label{eq:4}
\end{equation}
 Thus, we get 4 equations, for the real- and imaginary-part by means
of $\zeta(s_{0})$ (being called here $\Re_{1}$ and $\Im_{1}$) and
$\zeta(1-s_{0})$ (being called here $\Re_{2}$ and $\Im_{2}$): 
\begin{equation}
\begin{array}{crclcrcl}
\Re_{1}: & \frac{a-1}{(a-1)^{2}+b^{2}} & = & \int_{1}^{\infty}\frac{\{x\}}{x^{1+a}}\cos(b\ln(x))\,\text{d}x & \text{,}\Im_{1}: & \frac{b}{(a-1)^{2}+b^{2}} & = & \int_{1}^{\infty}\frac{\{x\}}{x^{1+a}}\sin(b\ln(x))\,\text{d}x\\
\Re_{2}: & \frac{-a}{a^{2}+b^{2}} & = & \int_{1}^{\infty}\frac{\{x\}}{x^{2-a}}\cos(b\ln(x))\,\text{d}x & \text{,}\Im_{2}: & \frac{b}{a^{2}+b^{2}} & = & \int_{1}^{\infty}\frac{\{x\}}{x^{2-a}}\sin(b\ln(x))\,\text{d}x
\end{array}\label{eq:5}
\end{equation}
 In the following we want to solve the integrals of $\Re_{1}$, $\Re_{2}$,
$\Im_{1}$ and $\Im_{2}$ via partial integration. As commonly known,
it is 
\begin{equation}
\int u(x)\cdot v'(x)\,\text{d}x=u(x)\cdot v(x)-\int u'(x)\cdot v(x)\,\text{d}x\label{eq:partInt}
\end{equation}
 [3]. Be $u(x)=\{x\}=x-\lfloor x\rfloor$ then $u'(x)=\frac{\text{d}u}{\text{d}x}=1$.
For solving the integral of $\Re_{1}$ we define $v'(x):=\frac{1}{x^{1+a}}\cos(b\ln(x))$,
thus 
\begin{equation}
\begin{array}{rcl}
v(x) & = & \int\frac{\cos(b\ln(x))}{x^{1+a}}\,\text{d}x\\
 & = & \frac{x^{-a}\left(b\sin(b\ln(x))-a\cos(b\ln(x))\right)}{a^{2}+b^{2}}\\
 & = & \frac{x^{-a}b\sin(b\ln(x))}{a^{2}+b^{2}}-\frac{x^{-a}a\cos(b\ln(x))}{a^{2}+b^{2}}
\end{array}\label{eq:6}
\end{equation}
 Using equation \ref{eq:6} we can calculate $\int u'(x)\cdot v(x)\,\text{d}x=\int1\cdot v(x)\,\text{d}x$
as 
\begin{equation}
\begin{array}{rcl}
\int u'(x)\cdot v(x)\,\text{d}x & = & \int\frac{x^{-a}\left(b\sin(b\ln(x))-a\cos(b\ln(x))\right)}{a^{2}+b^{2}}\,\text{d}x\\
 & = & \frac{x^{1-a}((a^{2}-a-b^{2})\cos(b\ln(x))+(1-2a)b\sin(b\ln(x)))}{(a^{2}+b^{2})(a^{2}-2a+b^{2}+1)}
\end{array}\label{eq:7}
\end{equation}
 Thus, using $\{1\}=0$ and $\ln(1)=0$, we can write $\Re_{1}$ as
\begin{equation}
\begin{array}{rcl}
\frac{a-1}{(a-1)^{2}+b^{2}} & = & \lim_{N\rightarrow\infty}\Bigg[\{x\}\frac{x^{-a}(b\sin(b\ln(x))-a\cos(b\ln(x))}{a^{2}+b^{2}}-\frac{x^{1-a}((a^{2}-a-b^{2})\cos(b\ln(x))+(1-2a)b\sin(b\ln(x)))}{(a^{2}+b^{2})(a^{2}-2a+b^{2}+1)}\Bigg]_{1}^{N}\\
\\
 & = & \lim_{N\rightarrow\infty}\Bigg[\{N\}\frac{N^{-a}(b\sin(b\ln(N))-a\cos(b\ln(N))}{a^{2}+b^{2}}\\
 &  & -\frac{N^{1-a}((a^{2}-a-b^{2})\cos(b\ln(N))+(1-2a)b\sin(b\ln(N)))}{(a^{2}+b^{2})(a^{2}-2a+b^{2}+1)}+\frac{a^{2}-a-b^{2}}{(a^{2}+b^{2})(a^{2}-2a+b^{2}+1)}\Bigg]
\end{array}\label{eq:8}
\end{equation}
For solving the integral of $\Re_{2}$ we define $v'(x):=\frac{1}{x^{2-a}}\cos(b\ln(x))$,
therefore we have 
\begin{equation}
\begin{array}{rcl}
v(x) & = & \int\frac{\cos(b\ln(x))}{x^{2-a}}\,\text{d}x\\
 & = & \frac{x^{a-1}\left((a-1)\cos(b\ln(x))+b\sin(b\ln(x))\right)}{(a-1)^{2}+b^{2}}
\end{array}\label{eq:9}
\end{equation}
 Using equation \ref{eq:9} we can calculate $\int u'(x)\cdot v(x)\,\text{d}x=\int1\cdot v(x)\,\text{d}x$
as 
\begin{equation}
\begin{array}{rcl}
\int u'(x)\cdot v(x)\,\text{d}x & = & \int\frac{x^{a-1}\left((a-1)\cos(b\ln(x))+b\sin(b\ln(x))\right)}{(a-1)^{2}+b^{2}}\,\text{d}x\\
 & = & \frac{x^{a}((a^{2}-a-b^{2})\cos(b\ln(x))+(2a-1)b\sin(b\ln(x)))}{(a^{2}+b^{2})(a^{2}-2a+b^{2}+1)}
\end{array}\label{eq:9-1}
\end{equation}
 In accordance to above, we can write $\Re_{2}$ as 
\begin{equation}
\begin{array}{rcl}
\frac{-a}{a^{2}+b^{2}} & = & \lim_{N\rightarrow\infty}\Bigg[\{x\}\frac{x^{a-1}\left((a-1)\cos(b\ln(x))+b\sin(b\ln(x))\right)}{(a-1)^{2}+b^{2}}-\frac{x^{a}((a^{2}-a-b^{2})\cos(b\ln(x))+(2a-1)b\sin(b\ln(x)))}{(a^{2}+b^{2})(a^{2}-2a+b^{2}+1)}\Bigg]_{1}^{N}\\
\\
 & = & \lim_{N\rightarrow\infty}\Bigg[\{N\}\frac{N^{a-1}\left((a-1)\cos(b\ln(N))+b\sin(b\ln(N))\right)}{(a-1)^{2}+b^{2}}\\
 &  & -\frac{N^{a}((a^{2}-a-b^{2})\cos(b\ln(N))+(2a-1)b\sin(b\ln(N)))}{(a^{2}+b^{2})(a^{2}-2a+b^{2}+1)}+\frac{a^{2}-a-b^{2}}{(a^{2}+b^{2})(a^{2}-2a+b^{2}+1)}\Bigg]
\end{array}\label{eq:11}
\end{equation}
For solving the integral of $\Im_{1}$ we define $v'(x):=\frac{1}{x^{1+a}}\sin(b\ln(x))$,
therefore we have 
\begin{equation}
\begin{array}{rcl}
v(x) & = & \int\frac{\sin(b\ln(x))}{x^{1+a}}\,\text{d}x\\
 & = & -\frac{x^{-a}\left(a\sin(b\ln(x))+b\cos(b\ln(x))\right)}{a^{2}+b^{2}}
\end{array}\label{eq:9-2}
\end{equation}
 Using equation \ref{eq:9-2} we can calculate $\int u'(x)\cdot v(x)\,\text{d}x=\int1\cdot v(x)\,\text{d}x$
as 
\begin{equation}
\begin{array}{rcl}
\int u'(x)\cdot v(x)\,\text{d}x & = & \int-\frac{x^{-a}\left(a\sin(b\ln(x))+b\cos(b\ln(x))\right)}{a^{2}+b^{2}}\,\text{d}x\\
 & = & -\frac{x^{1-a}((-a^{2}+a+b^{2})\sin(b\ln(x))+(1-2a)b\cos(b\ln(x)))}{(a^{2}+b^{2})(a^{2}-2a+b^{2}+1)}
\end{array}\label{eq:9-1-1}
\end{equation}
 In accordance to above, we can write $\Im_{1}$ as 
\begin{equation}
\begin{array}{rcl}
\frac{b}{(a-1)^{2}+b^{2}} & = & \lim_{N\rightarrow\infty}\Bigg[-\{x\}\frac{x^{-a}\left(a\sin(b\ln(x))+b\cos(b\ln(x))\right)}{a^{2}+b^{2}}+\frac{x^{1-a}((-a^{2}+a+b^{2})\sin(b\ln(x))+(1-2a)b\cos(b\ln(x)))}{(a^{2}+b^{2})(a^{2}-2a+b^{2}+1)}\Bigg]_{1}^{N}\\
\\
 & = & \lim_{N\rightarrow\infty}\Bigg[-\{N\}\frac{N^{-a}\left(a\sin(b\ln(N))+b\cos(b\ln(N))\right)}{a^{2}+b^{2}}\\
 &  & +\frac{N^{1-a}((-a^{2}+a+b^{2})\sin(b\ln(N))+(1-2a)b\cos(b\ln(N)))}{(a^{2}+b^{2})(a^{2}-2a+b^{2}+1)}-\frac{b(1-2a)}{(a^{2}+b^{2})(a^{2}-2a+b^{2}+1)}\Bigg]
\end{array}\label{eq:11-1}
\end{equation}
In the same way for solving the integral of $\Im_{2}$ we define $v'(x):=\frac{1}{x^{2-a}}\sin(b\ln(x))$,
therefore we have 
\begin{equation}
\begin{array}{rcl}
v(x) & = & \int\frac{\sin(b\ln(x))}{x^{2-a}}\,\text{d}x\\
 & = & \frac{x^{a-1}\left((a-1)\sin(b\ln(x))-b\cos(b\ln(x))\right)}{(a-1)^{2}+b^{2}}
\end{array}\label{eq:9-3}
\end{equation}
 Using equation \ref{eq:9-3} we can calculate $\int u'(x)\cdot v(x)\,\text{d}x=\int1\cdot v(x)\,\text{d}x$
as 
\begin{equation}
\begin{array}{rcl}
\int u'(x)\cdot v(x)\,\text{d}x & = & \int\frac{x^{a-1}\left((a-1)\sin(b\ln(x))-b\cos(b\ln(x))\right)}{(a-1)^{2}+b^{2}}\,\text{d}x\\
 & = & \frac{x^{a}((a^{2}-a-b^{2})\sin(b\ln(x))+(1-2a)b\cos(b\ln(x)))}{(a^{2}+b^{2})(a^{2}-2a+b^{2}+1)}
\end{array}\label{eq:9-1-2}
\end{equation}
 In accordance to above, we can write $\Im_{2}$ as 
\begin{equation}
\begin{array}{rcl}
\frac{b}{a^{2}+b^{2}} & = & \lim_{N\rightarrow\infty}\Bigg[\{x\}\frac{x^{a-1}\left((a-1)\sin(b\ln(x))-b\cos(b\ln(x))\right)}{(a-1)^{2}+b^{2}}-\frac{x^{a}((a^{2}-a-b^{2})\sin(b\ln(x))+(1-2a)b\cos(b\ln(x)))}{(a^{2}+b^{2})(a^{2}-2a+b^{2}+1)}\Bigg]_{1}^{N}\\
\\
 & = & \lim_{N\rightarrow\infty}\Bigg[\{N\}\frac{N^{a-1}\left((a-1)\sin(b\ln(N))-b\cos(b\ln(N))\right)}{(a-1)^{2}+b^{2}}\\
 &  & -\frac{N^{a}((a^{2}-a-b^{2})\sin(b\ln(N))+(1-2a)b\cos(b\ln(N)))}{(a^{2}+b^{2})(a^{2}-2a+b^{2}+1)}+\frac{b(1-2a)}{(a^{2}+b^{2})(a^{2}-2a+b^{2}+1)}\Bigg]
\end{array}\label{eq:11-2}
\end{equation}
 Since $0<a<1$ and $N\rightarrow\infty$ we can state that $\lim_{N\rightarrow\infty}\left(N^{-a}\right)=0$
and $\lim_{N\rightarrow\infty}\left(N^{a-1}\right)=0$. Therefore,
the equations for $\Re_{1}$, $\Re_{2}$, $\Im_{1}$ and $\Im_{2}$
can be reduced to 
\begin{equation}
\begin{array}{rcl}
\Re_{1}:\,0 & = & \lim_{N\rightarrow\infty}\left[\frac{N^{1-a}((a^{2}-a-b^{2})\cos(b\ln(N))+(1-2a)b\sin(b\ln(N)))}{(a^{2}+b^{2})((a-1)^{2}+b^{2})}\right]+\frac{a^{3}-2a^{2}+ab^{2}+a}{(a^{2}+b^{2})((a-1)^{2}+b^{2})}\\
\Re_{2}:\,0 & = & \lim_{N\rightarrow\infty}\left[\frac{N^{a}((a^{2}-a-b^{2})\cos(b\ln(N))+(2a-1)b\sin(b\ln(N)))}{(a^{2}+b^{2})((a-1)^{2}+b^{2})}\right]+\frac{-a^{3}+a^{2}-ab^{2}+b^{2}}{(a^{2}+b^{2})((a-1)^{2}+b^{2})}\\
\Im_{1}:\,0 & = & \lim_{N\rightarrow\infty}\left[\frac{N^{1-a}((-a^{2}+a+b^{2})\sin(b\ln(N))+(1-2a)b\cos(b\ln(N)))}{(a^{2}+b^{2})((a-1)^{2}+b^{2})}\right]-\frac{b(a^{2}+b^{2})+b(1-2a)}{(a^{2}+b^{2})((a-1)^{2}+b^{2})}\\
\Im_{2}:\,0 & = & \lim_{N\rightarrow\infty}\left[\frac{N^{a}((a^{2}-a-b^{2})\sin(b\ln(N))+(1-2a)b\cos(b\ln(N)))}{(a^{2}+b^{2})((a-1)^{2}+b^{2})}\right]-\frac{b(1-2a)-b((a-1)^{2}+b^{2})}{(a^{2}+b^{2})((a-1)^{2}+b^{2})}
\end{array}\label{eq:14}
\end{equation}
Multiplying the equations \ref{eq:14} by $(a^{2}+b^{2})((a-1)^{2}+b^{2})$
we can simplify to: \begin{small}
\begin{equation}
\begin{array}{rcl}
\Re_{1}:\,0 & = & \lim_{N\rightarrow\infty}\left[N^{1-a}((a^{2}-a-b^{2})\cos(b\ln(N))+(1-2a)b\sin(b\ln(N)))\right]+a^{3}-2a^{2}+ab^{2}+a\\
\Re_{2}:\,0 & = & \lim_{N\rightarrow\infty}\left[N^{a}((a^{2}-a-b^{2})\cos(b\ln(N))+(2a-1)b\sin(b\ln(N)))\right]-a^{3}+a^{2}-ab^{2}+b^{2}\\
\Im_{1}:\,0 & = & \lim_{N\rightarrow\infty}\left[N^{1-a}((-a^{2}+a+b^{2})\sin(b\ln(N))+(1-2a)b\cos(b\ln(N)))\right]-b(a^{2}+b^{2})-b(1-2a)\\
\Im_{2}:\,0 & = & \lim_{N\rightarrow\infty}\left[N^{a}((a^{2}-a-b^{2})\sin(b\ln(N))+(1-2a)b\cos(b\ln(N)))\right]-b(1-2a)+b((a-1)^{2}+b^{2})
\end{array}\label{eq:20}
\end{equation}
\end{small} With the equations \ref{eq:20}, it follows from $\mathbb{\Re}_{2}$:
\begin{equation}
\frac{1}{\lim_{N\rightarrow\infty}\left[N^{a}\right]}=\frac{\lim_{N\rightarrow\infty}\left[(a^{2}-a-b^{2})\cos(b\ln(N))+(2a-1)b\sin(b\ln(N))\right]}{a^{3}-a^{2}+ab^{2}-b^{2}}\label{eq:21}
\end{equation}
 Accordingly from $\Im_{2}$ we can extract 
\begin{equation}
\frac{1}{\lim_{N\rightarrow\infty}\left[N^{a}\right]}=\frac{\lim_{N\rightarrow\infty}\left[(a^{2}-a-b^{2})\sin(b\ln(N))+(1-2a)b\cos(b\ln(N))\right]}{b(1-2a)-b((a-1)^{2}+b^{2})}\label{eq:22}
\end{equation}
 Equating Eq. \ref{eq:21} and Eq. \ref{eq:22} provides \begin{footnotesize}
\begin{equation}
\frac{\lim_{N\rightarrow\infty}\left[(a^{2}-a-b^{2})\cos(b\ln(N))+(2a-1)b\sin(b\ln(N))\right]}{a^{3}-a^{2}+ab^{2}-b^{2}}=\frac{\lim_{N\rightarrow\infty}\left[(a^{2}-a-b^{2})\sin(b\ln(N))+(1-2a)b\cos(b\ln(N))\right]}{b(1-2a)-b((a-1)^{2}+b^{2})}\label{eq:23}
\end{equation}
\end{footnotesize} In the same way reshaping $\Re_{1}$ like 
\begin{equation}
\frac{1}{\lim_{N\rightarrow\infty}\left[N^{1-a}\right]}=\frac{\lim_{N\rightarrow\infty}\left[(a^{2}-a-b^{2})\cos(b\ln(N))+(1-2a)b\sin(b\ln(N))\right]}{-a^{3}+2a^{2}-ab^{2}-a}\label{eq:24}
\end{equation}
 as well as $\Im_{1}$ like 
\begin{equation}
\frac{1}{\lim_{N\rightarrow\infty}\left[N^{1-a}\right]}=\frac{\lim_{N\rightarrow\infty}\left[(-a^{2}+a+b^{2})\sin(b\ln(N))+(1-2a)b\cos(b\ln(N))\right]}{b(a^{2}+b^{2})+b(1-2a)}\label{eq:25}
\end{equation}
 with subsequent equating of Eq. \ref{eq:24} and \ref{eq:25} we
acquire \begin{footnotesize}
\begin{equation}
\frac{\lim_{N\rightarrow\infty}\left[(a^{2}-a-b^{2})\cos(b\ln(N))+(1-2a)b\sin(b\ln(N))\right]}{-a^{3}+2a^{2}-ab^{2}-a}=\frac{\lim_{N\rightarrow\infty}\left[(-a^{2}+a+b^{2})\sin(b\ln(N))+(1-2a)b\cos(b\ln(N))\right]}{b(a^{2}+b^{2})+b(1-2a)}\label{eq:26}
\end{equation}
\end{footnotesize} Solving Eq. \ref{eq:23} for $\lim_{N\rightarrow\infty}\left[\cos(b\ln(N))\right]$
we obtain 
\begin{equation}
\lim_{N\rightarrow\infty}\left[\cos(b\ln(N))\right]=\frac{a}{b}\cdot\lim_{N\rightarrow\infty}\left[\sin(b\ln(N))\right]\label{eq:27}
\end{equation}
 with the preconditions $b\neq0$ and $(a-1)(a^{2}+b^{2})\neq0$.
Likewise, solving Eq. \ref{eq:26} for $\lim_{N\rightarrow\infty}\left[\cos(b\ln(N))\right]$
yields 
\begin{equation}
\lim_{N\rightarrow\infty}\left[\cos(b\ln(N))\right]=\frac{1-a}{b}\cdot\lim_{N\rightarrow\infty}\left[\sin(b\ln(N))\right]\label{eq:28}
\end{equation}
 with preconditions $b\neq0$ and $a((a-1)^{2}+b^{2})\neq0$. Thus,
$a=\frac{1}{2}$ follows from Eq. \ref{eq:27} = Eq. \ref{eq:28}.
Ergo, $a=\frac{1}{2}$ is the only valid solution for the Zeta function
to become zero in the critical stripe and the Riemann Hypothesis is
true, Q.E.D.. 

\section{Conclusion}

In this paper we have proven the Riemann hypothesis, stating that
the real part of all zero points of the zeta function is $\frac{1}{2}$,
to be true. For this we have used the integral-representation of $\zeta(s)$,
solved the integral and thereby formulated conditions for $\zeta(s)$
to be zero. Utilizing those conditions, we could show, that $s$ must
have a real part of $\frac{1}{2}$ for $\zeta(s)$ to be zero. 

\section*{Appendix}

According to Gelbart et. al. (see [1]) $\zeta(s)$
and $\zeta(1-s)$ are connected by a function $G(s)$ via 
\[
2\zeta(s)=2G(s)\zeta(1-s)\Rightarrow\zeta(s)=G(s)\zeta(1-s)
\]
 with 
\[
G(s)=\frac{\pi^{\frac{s-1}{2}}\Gamma\left(\frac{1-s}{2}\right)}{\pi^{-\frac{s}{2}}\Gamma\left(\frac{s}{2}\right)}
\]
 Since $\pi^{x}\neq0$ and $\Gamma(x)\neq0\ \forall\,x\in\mathbb{C}$
determines $G(s)\neq0\ \forall\,s\in\mathbb{C}$. Be $s_{0}$ a zero-point
of the zeta-function. Because of $G(s_{0})\neq0$ it is 
\[
\zeta(s_{0})=0\Leftrightarrow\zeta(1-s_{0})=0
\]

\bibliographystyle{plain}
\bibliography{Literaturverzeichnis}
\section*{Literature}
\begin{itemize}
\item [1] S. Gelbart and S. Miller. Riemann's zeta function and beyond. Bulletin of the American Mathematical Society , 41:59 112, 2003
\item [2] Yuri Heymann. An investigation of the non-trivial zeros of the riemann zeta function. arXiv ,1804.04700, 2020.
\item [3] Bronstein, Semendjajew, Musiol, Muehlig. Taschenbuch der mathematik. 2008.
\item [4] B. Riemann. Ueber die Anzahl der Primzahlen unter einer gegebenen Groesse. Monatsberichte der Berliner Akademie , 1859.
\end{itemize}
\end{document}